\documentclass[12pt]{amsart}

\usepackage{amsmath}
\usepackage{amsthm}
\usepackage{amssymb}

\usepackage{beton,euler}

\usepackage{hyperref}

\author {Beno\^\i{}t Collins}
\address{Department of Mathematics and Statistics, 
University of Ottawa, 585 King Edward, K1N6N5
Ottawa, ON, Canada \&
CNRS, UMR 5208, Institut Camille Jordan, Universit\'e Lyon 1,
21 av Claude Bernard
69622  Villeurbanne, France}
\email{collins@math.univ-lyon1.fr}

\author{Piotr \'Sniady}
\address{Institute of Mathematics,
University of Wroclaw, pl.\ Grunwaldzki~2/4, 50-384
Wroclaw, Poland} \email{Piotr.Sniady@math.uni.wroc.pl}

\title[Representations of Lie groups and random matrices]
{Representations of Lie groups \\ and random matrices}

\theoremstyle{plain}
\newtheorem{lemma}{Lemma}
\newtheorem{theorem}[lemma]{Theorem}
\newtheorem{proposition}[lemma]{Proposition}
\newtheorem{corollary}[lemma]{Corollary}

\theoremstyle{definition}

\theoremstyle{remark}
\newtheorem{remark}[lemma]{Remark}
\newtheorem{example}[lemma]{Example}

\newcommand{\A}{{\mathfrak{A}}}

\newcommand{\E}{{\mathbb{E}}}
\newcommand{\C}{{\mathbb{C}}}
\newcommand{\R}{{\mathbb{R}}}
\newcommand{\Z}{{\mathbb{Z}}}
\newcommand{\g}{{\mathfrak{g}}}
\newcommand{\h}{{\mathfrak{h}}}
\newcommand{\uu}{{\mathfrak{u}}}

\newcommand{\El}{{\mathcal{L}}}
\newcommand{\M}{{\mathcal{M}}}

\newcommand{\gwia}{^{\star}}

\DeclareMathOperator{\tr}{tr} \DeclareMathOperator{\Tr}{Tr}

\DeclareMathOperator{\diag}{diag} 
\DeclareMathOperator{\End}{End} 
\DeclareMathOperator{\Id}{Id} 
\DeclareMathOperator{\Ad}{Ad}

\DeclareMathOperator{\Sym}{Sym}

\DeclareMathOperator{\Spin}{Spin}

\begin{document}

\begin{abstract}
We study the asymptotics of representations of a fixed compact Lie group.
We prove that the limit behavior of a sequence of such
representations can be described in terms of certain random matrices; in
particular operations on representations (for example: tensor product,
restriction to a subgroup) correspond to some natural operations on random
matrices (respectively: sum of independent random matrices, taking the corners
of a random matrix). Our method of
proof is to treat the canonical block matrix associated to a representation
as a random matrix with non-commutative entries.
\end{abstract}

\maketitle

\section{Introduction}

\subsection{Need for an asymptotic theory of representations}

One of the main questions in representation theory of Lie groups and Lie algebras  is to understand multiplicities, commutant spaces and structure
of representations arising in various natural situations, such as 
restriction to a subgroup or tensor product of representations.

There are many reasons to study asymptotic versions of such questions in the limit when the representation (and, possibly, also the Lie group) tends in some sense to infinity.
\begin{itemize}
\item  
From the viewpoint of probability theory it is
natural to consider the limit theorems (such as laws of large numbers, central limit theorem, etc.)
in order to study the limits 
of probability measures on a given set.
Reducible representations of a given group, the subject of this article,
can be alternatively described as probability measures on the set of
irreducible representations.

\item Even though for nearly all problems in the representation theory there are
explicit answers \cite{FultonHarris,GoodmanWallach}, they are
based on some combinatorial algorithms which are too cumbersome to be tractable
asymptotically. For this reason in the asymptotic theory of representations one
has to look for non-combinatorial tools such as random matrix theory
\cite{Mehta} or Voiculescu's free probability theory
\cite{VoiculescuDykemaNica}.

\item In theoretical physics it is a natural question as well. Indeed, finite
dimensional representation theory is described as a nice discrete object which
can be scaled in some thermodynamic limit to model continuous
phenomena \cite{KazakovStaudacherWynter1,KazakovStaudacherWynter2}.

\item Many important questions in the theory of operator algebras concern free
group factors. One of the foundations of Voiculescu's free probability theory
was the fact that free products may be approximated in some sense by random
matrices.  It was observed by Biane \cite{Biane95} that also representations
may provide such finite-dimensional approximants with an adventage of being
fully constructive and non-random.

\end{itemize}

\subsection{Asymptotics of representations of a fixed Lie group}
\label{subsec:asymptotics-representations-fixed}
Let $G$ be a fixed compact Lie group 
(in the following we consider the example when $G = U (d)$ is the unitary
group). 
One of our motivations is to 
study the problem of decomposing a
given concrete reducible representation of $G$ into irreducible components. 
For example we wish to study the following interesting examples of reducible
representations:
\begin{enumerate}
\item \label{problem-first} restriction of a given irreducible representation
of $G'$ to its subgroup $G$, where $G'$ is a given compact Lie group (for
example: $G = U (d)$ and $G' = U (d' )$ with $d < d'$);

\item \label{problem-last} Kronecker tensor product of given two irreducible
representations of $G$.
\end{enumerate}

The irreducible representations of compact Lie groups are uniquely determined
by their highest weights (cf Section
\ref{sec:representations} for definitions); in the example of $G = U (d )$
the irreducible representations are
indexed by sequences of integers $\lambda_1 \geq \cdots \geq
\lambda_d$. 
For explicit answers to the above questions 
there are well--known algorithms
involving combinatorial manipulations on the the highest weights. For example,
the decomposition of an irreducible representation of $G' = U ( d+ 1 )$ to a
subgroup $G = U (d )$ is given by the Weyl's branching rule; in the general
case the answer to this question is given by counting certain Littelmann paths
\cite{Littelmann95}.
However, when we are interested in the situation when the dimension of the
representations tends to infinity, these combinatorial algorithms become very
cumbersome. In
particular, the direct study of the multiplicities in the above problems seems
to be rather difficult.

In order to avoid such difficulties we 
concentrate on approximate asymptotic answers.
More explicitly, for a given sequence $(\rho_n)$
of representations of $G$ for which the highest weights of the irreducible
components tend to infinity, we study the limit distribution of a
(rescaled) highest weight of a randomly chosen irreducible component of
$\rho_n$.

This situation in the asymptotic theory of representions in many ways
resembles statistical mechanics: when the number of particles in a
physical system grows, the complexity of its description also grows so that its
exact solution becomes quickly intractable.  
However, a more modest approach in
which we are interested only in some macroscopic quantities may result in a
surprisingly simple description. Similarly, in the asymptotic theory of
representations when we abandon the attempts to find exact combinatorial
solutions and restrict ourselves to a statistical description we may expect
dramatic simplifications. Since representations have a highly non-commutative
structure and the combinatorial algorithms behind are so cumbersome therefore we
can expect quite surprising results.

\subsection{The main result: Representations of Lie groups and random matrices}

For a fixed compact Lie group $G$ and a sequence $(\rho_n)$ of 
representations we construct a random matrix whose joint eigenvalues
distribution depends on the asymptotic behavior of the
highest weights of $\rho_n$. The symmetry of this random matrix depends on
the group $G$; for example when $G=U(d)$ is the unitary group the corresponding
random matrix will be a hermitian $d\times d$ random matrix which additionally
is invariant under conjugation by unitary matrices.

We prove that some operations on the representations (such as
Kronecker tensor product, restriction to a subgroup) correspond to some natural
operations on the appropriate random matrices (sum of independent random
matrices, taking the corner of a random matrix, respectively). In this way
problems concerning asymptotic properties of representations are reduced to
much simpler analytic problems concerning random matrices.

Some results proved in this paper were already 
 considered by Heckman \cite{Heckman1982}. However, as far as we could understand,
 his proofs are very different from ours. Our methods are probabilistic and lead to many 
 new applications and examples, see Section \ref{sec:applications}.


\subsection{The main tool: random matrices with non-commutative entries}
Let $\rho$ be a representation of a compact Lie group $G$; from the
following on we shall restrict our attention to the corresponding representation
of the Lie algebra $\g$ (if $G$ is connected the
representation of the Lie group is uniquely determined by the representation of
the Lie algebra).

The family of matrices $\{ \rho(x): x\in \g\}$ can be
viewed as a family of \emph{non-commutative random variables}; in other words
$\rho$ can be viewed as a \emph{non-commutative random vector} in $\g\gwia$.
We prove that asymptotically, when the representation $\rho$ tends to
infinity, such a non-commutative random vector converges in distribution (after
some rescaling) to a classical (\emph{commutative}) random vector in $\g\gwia$
and hence---in many cases---can be identified with a random matrix.

This idea is closely related to the work of Kuperberg
\cite{Kuperberg2002,Kuperberg2005} who---among other results---gives a new,
conceptual proof of the result of Johansson \cite{Johansson01} (see Theorem
\ref{theo:clt} below). Kuperberg's idea is to treat elements of the Lie
algebra as non-commutative random variables and to show that for the tensor
product of representations a non-commutative central limit theorem can be
applied. The results of this article can be therefore viewed as an extension of
some of Kuperberg's results \cite{Kuperberg2002} from central limit theorem
related to the tensor product of representations to other operations on
representations such as restrictions or tensor products of a fixed number of
representations.

\subsection{Asymptotics of representations for a series of Lie groups and free
probability}
\label{subsec:asympt-series}

A variation of the above problem appears when we replace the fixed group $G$ by a classical series of compact Lie groups $(G_n)$ and we consider a series $(\rho_n)$, where $\rho_n$ is a representation of $G_n$; we are interested in solving the analogues of the problems
(\ref{problem-first})--(\ref{problem-last}) from Section \ref{subsec:asymptotics-representations-fixed}. For example, we may take $G_n=U(n)$ to be the series of the unitary groups. This case was studied in detail by Biane \cite{Biane95} who found a connection between asymptotics of such representations and Voiculescu's free probability theory \cite{VoiculescuDykemaNica}. 

In a subsequent paper \cite{CollinsSniady2007} we show that the method of random matrices with non-commuting entries can also be applied to this situation and the results obtained in this way are significantly stronger than the ones of Biane \cite{Biane95}.

\subsection{Asymptotics of representations of symmetric groups}
It turns out that the ideas presented in this article can be also applied to the case of the symmetric groups. A canonical matrix associated to a representation of the symmetric group was given by Biane \cite{Biane1998} and it turns out that the recent results of the second-named author \cite{Sniady2005GaussuanFluctuationsofYoungdiagrams} were proved by treating (in a very concealed way) this matrix as a permutationally-invariant random matrix with non-commuting entries. A subsequent paper \cite{SniadySpeicher2007permutationmatrix} will present the details.

\subsection{Overview of this article}

In Section \ref{sec:non-commutative} we introduce some notations concerning non-commutative random variables. In Section \ref{sec:representations} we state recall some facts about
representation theory and fix some notation. In Section \ref{sec:rep} we prove the main result and in Section \ref{sec:applications} we investigate its new consequences.

\section{Non-commutative random variables}
\label{sec:non-commutative}

\subsection{Non-commutative probability spaces}
Let $( \Omega, \mathfrak{M}, P )$ be a Kolmogorov probability space. We
consider an algebra
\[ \El^{\infty -} ( \Omega ) = \bigcap_{n\geq 1} \El^n(\Omega)\]
of random variables with all moments finite. This algebra is equipped with a
functional $\E : \El^{\infty -} ( \Omega ) \rightarrow \R$ which to a random
variable associates its mean value.

We consider a generalization of the above setup in which the commutative algebra
$\El^{\infty -} ( \Omega )$ is replaced by any (possibly non-com\-mu\-ta\-ti\-ve)
$\star$-algebra $\A$ with a unit and $\E : \A \rightarrow \R$ is any linear
functional which is normalized (i.e.\ $\E ( 1 ) = 1$), positive (i.e.\ $\E(x\gwia x)>0$ for all $x\in\A$ such that $x\neq 0$) and tracial (i.e.\ $\E(xy)=\E(yx)$ for all
$x,y\in\A$). 
The elements of $\A$
are called \textit{non-commutative random variables} and the functional $\E$
is called mean value or expectation. We also say that $( \A, \E )$ is
a non-commutative probability space \cite{VoiculescuDykemaNica}.

Here are two motivating examples which will be used in the following.
\begin{example}
  \label{ex:commutative}
For any Kolmogorov probability space the corresponding pair $\big( \El^{\infty
-} ( \Omega ), \E \big)$ is a non-commutative probability space.
\end{example}

\begin{example}
\label{ex:finite-matrices}
For any Kolmogorov probability space and integer $d\geq 1$ we consider the
algebra $ \El^{\infty -}(\Omega) \otimes \M_d = \El^{\infty -}(\Omega; \M_d )$
of $d\times d$ random matrices and we equip it with a functional
$$ \E(x) = E(\tr x), $$
for any random matrix $x$, where $E$ on the
right-hand-side denotes the mean value and
$$ \tr x = \frac{1}{d} \Tr x $$
is the normalized trace. In this way
$\big(\El^{\infty-}(\Omega; \M_d ),\E\big)$ is a non-commutative probability
space.
\end{example}

Any non-commutative probability space $(\A,\E)$ can be equipped with the
corresponding $L^2$ norm:
$$ \| x \|_{L^2} = \sqrt{\E ( x x\gwia )}. $$

Notice that the above definitions of non-commutative probability spaces and
random variables do not require any analytic notions other than positivity. 
In particular, as we shall see in the remaining part of this section, by the
distribution of a non-commutative random variable we understand the collection
of its mixed moments. While this approach turns out to to be very useful to
state and prove our results and it has the advantage to encompass many
non-bounded random variables, it has a couple of drawbacks. For instance, the
convergence of non-commutative distributions as defined in Section
\ref{sec:randomvectors} does not coincide in the commutative case with the weak
convergence of probability measures, therefore some of our corollaries
concerning convergence in distribution of classical random variables will be
formulated and proved not in the most desirable weak topology of probability
measures but with respect to the moments convergence. This issue and the way to
fix it in the cases which are of our interest (so that the convergence in the
weak topology of probability measures in fact holds true) are discussed in
Section \ref{sec:wtuniqueness}.

\subsection{Random vectors}
\label{sec:randomvectors}

Let $V$ be a finite dimensional (real) vector space. If
$v : \Omega \rightarrow V$ is a random variable valued in the space $V$ we
say that \textit{$v$ is a (classical) random vector in $V$}. We say that
$v$ has all moments finite if $\E \| v \|^k < \infty$ holds
true for any exponent $k \geq 1$. Notice that this definition does not depend
on the choice of the norm $\| \cdot \|$ on $V$.
For a random vector $v$ with all moments finite we define its \textit{moments}
\begin{equation}
  \label{eq:momenty-klasyczne} m_k = m_k^{\E} ( v ) = \E \underbrace{v \otimes
  \cdots \otimes v}_{k \text{ times}} \in V^{\otimes k} .
\end{equation}
In the case when $V = \R$, $v$ becomes a usual
number-valued random variable; furthermore $V^{\otimes k} = \R^{\otimes k}
\cong \R$ and the moments $m_k = \E v^k \in \R$ are just real numbers and this
definition coincides with the usual notion of the moments of a random variable.
In the following we are interested in the space
\[ \left\{ v : \Omega \rightarrow V \text{ such that } \E \| v \|^k < \infty
   \text{ for each } k \geq 1 \right\} \]
of random vectors with all moments finite which we will view as a tensor
product
\begin{equation}
  \label{eq:random-vector} V \otimes \El^{\infty -} ( \Omega ) .
\end{equation}

\subsection{Non-commutative random vectors}

Let $( \A, \E )$ be a non-com\-mu\-ta\-ti\-ve probability space; in analogy to
{\eqref{eq:random-vector}} we call the elements of $V \otimes \A$
\textit{non-com\-mu\-ta\-ti\-ve random vectors in $V$ (over a non-commutative
probability space $( \A, \E )$)}.

Given $v_1=x_1\otimes a_1\in V_1\otimes \A$ and $v_2=x_2\otimes a_2 \in
V_2\otimes \A$ we define 
$$v_1\hat{\otimes} v_2=(x_1\otimes a_1)\hat{\otimes} (x_2\otimes a_2)=
(x_1\otimes x_2\otimes a_1a_2)\in V_1\otimes V_2\otimes \A
$$
and its linear extension on non-elementary tensors.
Whenever $v_1=v_2$ with $V_1=V_2$ one shortens the notation 
as $v^{\hat{\otimes} 2}\in V^{\otimes 2}\otimes \A$ and one extends it 
by recursion to the definition of
$$v^{\hat{\otimes} k}\in V^{\otimes k}\otimes \A .$$
Observe that this definition matches the definition of the tensor product
of compact quantum groups of Woronowicz
\cite{Woronowicz87Compact-matrix-pseudogroups}
provided that $\A$ is a quantum group and $V$ a representation
of $\A$.

The \textit{$k$-th order vector moment} $m_k^{\E}(v)$ is defined as
$$m_k^{\E}(v)=(\Id\otimes  \E )v^{\hat{\otimes} k}\in 
V^{\otimes k}.$$
If there is no ambiguity we might remove the superscript $\E$.
Note that in the case of Example \ref{ex:commutative} when
$\A$ is commutative the moments, as defined above, coincide with the old definition
{\eqref{eq:momenty-klasyczne}} of the moments of a random vector.

We define the distribution of a non-commutative random vector as its sequence
$(m_k^{\E}(v))_{k=1,2,\dots}$ of moments. Accordingly, convergence in distribution of
non-commutative random vectors is to be understood as convergence of the
moments.

The above definitions can be made more explicit as follows: let $e_1,\dots,e_d$
be a base of the vector space $V$. Then a (classical) random vector $v$ in $V$
can be viewed as
\begin{equation}
\label{eq:wspolrzedne}
 v= \sum_i a_i e_i,
\end{equation}
where $a_i$ are the (random) coordinates. Then a non-commutative random vector
can be viewed as the sum \eqref{eq:wspolrzedne} in which $a_i$ are replaced by
non-commutative random variables. One can easily see that the sequence of
moments
$$ m_k(v) = \sum_{i_1,\dots,i_k} \E(a_{i_1} \cdots a_{i_k}) \    e_{i_1} \otimes
\cdots \otimes e_{i_k}$$
contains nothing else but the information about the mixed moments of the
non-commutative coordinates $a_1,\dots,a_d$ and the convergence of moments is
equivalent to the convergence of the mixed moments of $a_1,\dots,a_d$.

The following result provides a necessary and sufficient condition for a sequence of moments
to be those of a commutative vector.

\begin{proposition}
A non-commutative random vector $v$ actually arises from a commutative probability space 
iff for each value of $k\in\{1,2,\dots\}$ the tensor  $m_k^{\E}(v)\in V^{\otimes k}$ is 
invariant under the action of the symmetric group.
\end{proposition}
\begin{proof}
The necessity is trivial.
For sufficiency, if $m_k^{\E}(v)$ are invariant under the action by conjugation
of the symmetric group, this implies that the GNS representation of $\A$ with
respect to $\E$ is its abelianized quotient. Since $\E$ is supposed to be
faithful the proof is complete.
\end{proof}

\section{Preliminaries of representation theory} \label{sec:representations}

\subsection{Structure of compact Lie groups}

In this section we recall some facts about
Lie groups and their algebras \cite{Brocker-tomDieck,FultonHarris}.
Let $G$ be a compact Lie group and let $H \subseteq G$ be a
maximal abelian subgroup and let $\h \subseteq \g$ be the corresponding Lie
algebras. We consider the adjoint action of $\h$ on the vector space $\g$;
then there is a decomposition
\begin{equation}
\label{eq:decomposition}
\g = \h \oplus \bigoplus_{\alpha \in \h \gwia} \g_{\alpha}
\end{equation}
into eigenspaces. The non-zero elements $\alpha$ for which the
corresponding eigenspace $\g_{\alpha}$ is non-trivial are called \emph{roots}.
We identify $\h \gwia$ as a set of functionals on $\g$ which vanish on all root
spaces $\g_{\alpha}$ thus
\begin{equation}
\label{eq:inclusion}
\h \gwia \subseteq \g \gwia.
\end{equation}

Suppose that $\g$ is equipped with a
$G$-invariant scalar product; in this way $\g \cong \g\gwia$ and $\h \cong
\h\gwia$. The above isomorphisms and the inclusion $\h\subseteq\g$ allow us to
consider the inclusion $\h\gwia\subseteq \g\gwia$. This
inclusion does not depend on the choice of the $G$-invariant scalar product on $\g$ and it
coincides with the inclusion \eqref{eq:inclusion} considered above.
More generally, throughout all this paper, whenever we state a result related to
$\g$ but do  not mention the invariant scalar product, it means that the result
does  not depend on the choice
of the scalar product.

By definition, the Weyl group $W=\{ g\in G: g H g^{-1} = H \}$ is the set of the elements which normalize $H$. The Weyl group $W$ is always finite.

Every element $x\in\g$ is conjugate to some element $y\in\h$ (i.e.\ there
exists an element $g\in G$ such that $\Ad_g(x)=y$, where
$\Ad:G\rightarrow\End(\g)$ denotes the adjoint representation) and $y$ is
determined only up to the action of the Weyl group. Similarly, every element
$x\in \g \gwia$ is conjugate to some element $y\in\h \gwia$ (which is
determined only up to the action of the Weyl group).

\begin{example}
\label{example:Un}
For $G=U(n)$ we may take $H$ to be the group of diagonal unitary matrices.
Then $\g$ is the set of $n\times n$ antihermitian matrices and $\h$ is
the set of diagonal matrices with imaginary entries. In this case $W$ is the
group of permutation matrices and hence is isomorphic to the symmetric group
$S_n$. We equip $\g$ with a $G$-invariant scalar product $\langle
x,y\rangle=\Tr x y\gwia$ and thus we identify $\g\cong \g\gwia$ and $\h\cong
\h\gwia$. Clearly, every element  $x\in \g\cong \g\gwia$ is conjugate by a
unitary matrix to some diagonal matrix $y\in\h\cong
\h\gwia$ which is determined by the eigenvalues of $x$. The freedom of choosing
the order of the diagonal elements of $y$ corresponds to the action of the Weyl
group $S_n$.
\end{example}

\begin{remark}
\label{re:Un}
In random matrix theory it is more customary to work with the space of hermitian matrices instead of the space of antihermitian matrices; for this reason we may consider a simple isomorphism between these spaces given by multiplication by $i$. 
\end{remark}

The above example heuristically motivates that we may think that for a given
element $x\in \g\gwia$ the corresponding element $y\in \h\gwia$ contains the
information about the \emph{eigenvalues} of $x$.

\subsection{Irreducible representations and highest weights}

Let $\rho:\g\rightarrow \End(V)$ be an
irreducible representation. 
The representation space can be decomposed into the eigenspaces of the maximal
abelian subalgebra $\h$:
\begin{equation}
\label{eq:decomposition2}
 V = \bigoplus_{\alpha\in \h\gwia} V_{\alpha},
\end{equation}
in other words $h\in\h$ acts on $V_{\alpha}$ as multiplication by $\alpha(h)$.
The elements $\alpha\neq 0$ for which the corresponding eigenspace $V_{\alpha}$
is non-trivial are called \emph{weights}. 

Let us fix some element $t\in \h$ which
is generic in a sense that for any weight non-zero $\alpha$ we have
$\alpha(t)\neq 0$. Now a weight $\alpha$ is called \emph{positive} if
$\alpha(t)>0$, otherwise it is called \emph{negative}; similar convention
concerns roots as well. The \emph{highest weight} of a representation is the
weight $\alpha$ for which $\alpha(t)$ takes the maximal value. In fact, any
irreducible representation is uniquely determined (up to equivalence) by its
highest weight. The weight space $\g_{\alpha}$ corresponding to the highest
weight is always one-dimensional. Clearly, the above definitions depend on the choice of $t$, nevertheless the notion of highest weight is well-defined (up to the action of the Weyl group).

\subsection{Enveloping algebra}
The enveloping algebra of $\g$, denoted by $\bar{\g}$, is the free algebra
generated by the elements of $\g$ quotiented by relations $gh-hg=[g,h]$ for any
$g,h\in\g$. Usually we work over the complex (complexified) algebra. This
algebra is naturally endowed with a $\star$-algebra structure obtained by declaring
that the elements of $\g$ (before complexification) are anti-hermitian elements.
For our purposes, it has the  following
two important features:
\begin{itemize}
\item an irreducible representation $(\rho ,V)$ of $\g$ gives rise to
an onto algebra homomorphism $\bar{\g}\to \End(V)$ and conversely;
\item therefore, irreducible representations of $\g$ are in one to one
correspondence
with minimal tracial states on $\bar{\g}$.
\end{itemize}

We will need the following theorem, known as
Poincare-Birkhoff-Witt theorem.
\begin{theorem}
If\/ $h_1,\ldots h_n$ is a basis of\/ $\g$
as a vector space, then
$$(h_1^{\alpha_1}\cdots h_n^{\alpha_n})_{\alpha_1,\ldots, \alpha_n \in\mathbb{N}}$$
is a basis of\/ $\bar{\g}$ as a vector space. In particular, there is a
filtration on $\bar{\g}$ defined as follows:
the degree of  $p\in\bar{g}$ is the smallest $k$ such
that $p$ is a sum of monomials of elements in $g$ with at most $k$ factors.
\end{theorem}

\subsection{Reducible representations and random highest weights}

Let $\rho:G\rightarrow \End(V)$ be a (possibly reducible) representation of $G$
on a
finite-dimensional vector space $V$. We may decompose $\rho$ as a sum of
irreducible representations:
$$ \rho= \bigoplus_{\lambda\in \h\gwia} n_{\lambda} \rho_{\lambda}, $$
where $\rho_{\lambda}$ denotes the irreducible representation of $G$ with the
highest weight $\lambda$ and $n_{\lambda}\in\{0,1,2,\dots\}$ denotes its
multiplicity in $\rho$.

We define a probability measure $\mu_{\rho}$ on $\h\gwia$ such that the
probability of $\lambda\in\h\gwia$ is equal to
$$ \frac{n_{\lambda} \cdot (\text{dimension of
$\rho_{\lambda}$})}{\text{(dimension of $V$)}}; $$
in other words it is proportional to the total dimension of all the
summands of type $[\rho_{\lambda}]$ in $\rho$. In this way the probability
measure $\mu_{\rho}$ encodes in a compact way the information about the
decomposition of $\rho$ into irreducible components.
We define the \emph{random highest weight} associated to representation $\rho$
as a random variable distributed accordingly.

\begin{remark}
If we literally follow the above definition then $\mu_{\rho}$ 
is a certain probability measure on the Weyl chamber. This definition
has a disadvantage that it depends on the choice of the Weyl chamber, therefore
sometimes it will be convenient to understand by $\mu_{\rho}$ 
the $W$-invariant probability measure on $\h\gwia$ obtained by symmetrizing the above
measure by the action of the Weyl group $W$.
\end{remark}

\subsection{Random matrix with specified eigenvalues}
For an element $\lambda \in \h \gwia$ we consider a random vector in $\g
\gwia$ given by $\tilde{\lambda} = \Ad_g \lambda \in \g \gwia$, where $g$ is a
random element of $G$, distributed according to the Haar measure on $G$.
We will say that $\tilde{\lambda}$ is a \emph{$G$-invariant random matrix with
the eigenvalues $\lambda$}. This terminology was motivated by the following
example.

\begin{example}
If $G=U(n)$ is the group of the unitary matrices and
$\lambda=\diag(\lambda_1,\dots,\lambda_n)\in\h\gwia$ is a diagonal matrix then
the distribution of $\tilde{\lambda}$ is indeed equal to the uniform measure on
the manifold of antihermitian matrices with the eigenvalues
$\lambda_1,\dots,\lambda_n$.
\end{example}

We may extend the above definition to the case when $\lambda\in \h\gwia$ is
random; in such a case we will additionally assume that $g$ and
$\lambda$ are independent. If $\mu$ is a probability measure on $\h\gwia$  we
may treat it as distribution of the random variable $\lambda$; in this case we
we say that
$\tilde{\lambda}$ is a \emph{$G$-invariant random matrix with
the distribution of eigenvalues given by $\mu$}. We denote by
$\tilde{\mu}$ the corresponding distribution of $\tilde{\lambda}$.

\begin{remark}
\label{rem:sometimes_it_is_better_to_symmetrize}
It would be more appropriate to call $\lambda$ a $G$-invariant random
\emph{vector} in $\g\gwia$, nevertheless in the most interesting examples the
Lie groups under consideration carry some canonical matrix structure hence
the elements of $\g\gwia$ can be indeed viewed as matrices.
\end{remark}

The moments of $G$-invariant matrices are characterized by the following lemma.

\begin{lemma}
\label{lem:characterization}
Let $\mu$ be a probability measure on $\h \gwia$ with all moments finite and
let $\tilde{\mu}$ be the corresponding distribution of a $G$-invariant random
matrix.

For each $k$ the moment $m_k = m_k ( \tilde{\mu} ) \in (\g \gwia)^{\otimes k}$
is the unique element such that:
\begin{enumerate}

  \item \label{prop2} $m_k$ is invariant under the adjoint action of\/ $G$,

  \item \label{prop3}actions of\/ $m_k: \g^{\otimes k}
  \rightarrow \R$ and $m_k(\mu): \g^{\otimes k}
  \rightarrow \R$ coincide on $G$-invariant tensors in $\g^{\otimes k}$ (where
above $m_k(\mu)$ denotes the canonical extension of\/
$m_k(\mu)\in(\h\gwia)^{\otimes k}$ which is possible thanks to the inclusion
$\h\gwia\subseteq \g\gwia$).
\end{enumerate}
\end{lemma}

\begin{proof}
Point (\ref{prop2}) follows easily from the invariance of the Haar measure.

For any elements $x_1,x_2\in\g\gwia$ which are conjugate to each other the
restrictions of the maps $x_1^{\otimes k},x_2^{\otimes k}:\g^{\otimes
k}\rightarrow\R$ to $G$-invariant tensors coincide.
In this way point (\ref{prop3}) follows easily.

For $x\in \g^{\otimes k}$ let 
$$x'= \int_G \Ad_g^{\otimes k}(x) \ dg$$
be the average over the Haar measure on $G$.  Clearly, $x'$ is a $G$-invariant tensor.
Since
$$\big(m_k(\tilde{\mu}) \big)(x) = \big(m_k(\tilde{\mu}) \big)(x') $$
therefore the values of the functional $m_k(\tilde{\mu})$ are uniquely
determined by its values on $G$-invariant tensors which shows the uniqueness of
$m_k$.
\end{proof}

\begin{remark}
Similarly as in Remark \ref{rem:sometimes_it_is_better_to_symmetrize},
it is sometimes more convenient to understand the distribution of the eigenvalues of 
a random element of $\g\gwia$ as a $W$-invariant measure on $\h\gwia$.
\end{remark}

\section{Representations and random matrices with non-commutative entries}
\label{sec:rep}

\subsection{The main result}

Let $\rho : \g \rightarrow \End ( V )$ be a representation of $\g$.
We shall view $\rho$ as an element of $\g \gwia \otimes \End ( V )$; in other
words $\rho$ is a non-commutative random vector in $\g \gwia$ over a
non-commutative probability space $\big( \End (V), \tr \big)$.
The sequence of its
moments $m_k= m_k ( \rho ) \in ( \g \gwia )^{\otimes k}$, or equivalently, $m_k
: \g^{\otimes k} \rightarrow \C$ is given explicitly by
$$m_k( g_1 \otimes \cdots \otimes
g_k ) = \tr \big[ \rho ( g_1 ) \cdots \rho ( g_k
   ) \big].$$

The following theorem is the main result of this article.

\begin{theorem}
  \label{theo:main2}
Let $( \epsilon_n )$ be a sequence of real numbers which
converges to zero.  For each $n$ let $\rho_n : \g \rightarrow \End(V_n)$ be a
representation of $\g$ and let $\lambda_n$ be the corresponding random
highest weight. Let $A$ be a
$G$-invariant random matrix with all moments
finite.

Then the following conditions are equivalent:

\begin{enumerate}
\item \label{convergence:A} 
the distributions of the random variables
$\epsilon_n \lambda_n$ converge in moments to the distribution of the
eigenvalues of $A$;

\item \label{convergence:B} 
the sequence $\epsilon_n \rho_n$ of
non-commutative random matrices converges in distribution to $A$.
\end{enumerate}
\end{theorem}

\begin{proof}

Suppose that condition \eqref{convergence:A} holds true. It is enough to prove that from any
subsequence $(\epsilon_{k(n)} \rho_{k(n)})$
one can choose a subsequence $(\epsilon_{k(l(n))} \rho_{k(l(n))})$ which
converges in distribution to the random matrix $A$.

Let a subsequence $(\epsilon_{k(n)} \rho_{k(n)})$ be given; by Lemma
\ref{lem:operator-norm} and the compactness argument it follows that there
exists a subsequence $(\epsilon_{k(l(n))} \rho_{k(l(n))})$ which converges in
distribution to some random matrix $M \in \g\gwia \otimes \A$ with
non-commutative entries (for some non-commutative probability space $(\A,\E)$).
It is enough to prove that $M$ is the $G$-invariant random matrix in $\g \gwia$
with the same distribution of the eigenvalues as for $A$. In order to keep the
notation simple instead of $(\epsilon_{k(l(n))} \rho_{k(l(n))})$ we will write
$(\epsilon_{n} \rho_{n})$.

Firstly, observe that Lemma \ref{lem:operator-norm} shows that
\begin{multline*} \big\| [ \epsilon_n \rho_n(x_1), \epsilon_n \rho_n(x_2)] \big\|_{L^2}=
\big\| \epsilon^2_n \rho_n\big( [x_1,x_2] \big) \big\|_{L_2} \leq \\
   \big( \E \|\epsilon_n \lambda_n\|^2 \big)^{\frac{1}{2}} \ \epsilon_n 
\big\|[x_1,x_2]\big\| \to 0 
\end{multline*}
where we used the fact that the first factor on the right-hand side converges to
some constant depending on the distribution of eigenvalues of $A$. This shows that the 
elements $\{ M(x)\in\A : x\in\g \}$ commute hence $M$ can be identified with a
classical random variable (valued in $\g$).

Secondly, since each of the random matrices $\epsilon_n \rho_n$ is
$G$-invariant,
the random matrix $M$ is also $G$-invariant.

Thirdly, let $x\in\g^{\otimes k}$ be invariant under the action of
$\Ad^{\otimes k}:G\rightarrow\End(\g^{\otimes k})$. For every irreducible
representation $\rho$ we have that $\rho^{\otimes k}(x)$ is a multiple of
identity hence can be identified with a complex number. The exact value of this
number is equal to $\rho^{\otimes k}(x)\big|_{V_{\lambda}}$ which can be
estimated with the help of Lemma \ref{lem:highest-weight}. Therefore for a
(possibly reducible) representation $\rho_n$ we have
$$ \big( m_k(\epsilon_n \rho_n) \big)(x) = 
\E \Lambda_n^{\otimes k}(x) + (\text{terms of degree at least $1$ in
$\epsilon$}), $$
where $\Lambda_n=\epsilon_n \lambda_n$, where $\lambda_n$ is the random highest
weight associated to the representation $\rho_n$. Hence 
\begin{multline*} \big(m_k(M)\big)(x)= \lim_{n\to\infty }\big( m_k(\epsilon_n \rho_n) \big)(x)
= \\
\lim_{n\to\infty} \big( m_k (\epsilon_n \lambda_n) \big) (x) = \big( m_k
(A) \big) (x). 
\end{multline*}
The above equality and Lemma \ref{lem:characterization} show that the distribution of the eigenvalues of $M$
coincides with the distribution of eigenvalues of $A$ which finishes the proof.

Suppose that condition (\ref{convergence:B}) holds true. 
In order to prove (\ref{convergence:A})---similarly as in the proof
of the implication (\ref{convergence:A})$\implies$(\ref{convergence:B})---it is 
enough to show that from any subsequence $(\epsilon_{k(n)} \lambda_{k(n)})$
one can choose a subsequence $(\epsilon_{k(l(n))} \lambda_{k(l(n))})$ which
converges in moments to the distribution of the eigenvalues of $A$. 

Lemma \ref{lem:oszacowanie-lambdy} can be used to show that if $k$ is even 
then the sequence of moments 
$m_k(\epsilon_n \lambda_n)$ is bounded hence for every $k$
the sequence of moments $m_k(\epsilon_n \lambda_n)$ is bounded. 
Again, by a compactness argument we can find a subsequence which converges in moments to 
the distribution of eigenvalues of some random matrix $A'$; from the implication 
(\ref{convergence:A})$\implies$(\ref{convergence:B}) 
it follows that random matrices $A$ and $A'$ must have equal moments of their entries. 
From Lemma \ref{lem:entries_determine_eigenvalues} it follows that the eigenvalues of
$A'$ have the same moments as the eigenvalues of $A$ which finishes the proof.
\end{proof}

\subsection{Key lemmas}
\label{sec:proof}

We start with the following estimate:

\begin{lemma}
\label{lem:operator-norm}
We equip $\g$ with a $G$-invariant scalar product; in this way
we equip $\g\cong\g\gwia$ with the corresponding norm.

Let a unitary representation of a group $G$ on a finite-di\-men\-sio\-nal Hilbert
space be given and let $\rho$ be the corresponding representation of a Lie
algebra $\g$. If $\rho$ is irreducible with highest weight $\lambda$ then
for any $x\in\g$
\begin{equation}
\label{eq:operator-norm}
 \| \rho(x) \| \leq  \|x\| \ \|\lambda\|,
\end{equation}
where the norm on the left-hand side denotes the operator norm.
\end{lemma}
\begin{proof}
Since $x\in\g$ is conjugate to some element of $\h$, it is enough to prove
\eqref{eq:operator-norm} for $x\in\h$. For such elements the action of
$\rho(x)$ is diagonal with respect to the decomposition \eqref{eq:decomposition2}
and
$$ \| \rho(x) \| = \max_{\alpha} |\alpha(x)| \leq |\lambda(x)|, $$
where the maximum runs over the set of roots $\alpha\in\g\gwia$ contributing to
\eqref{eq:decomposition2}.
\end{proof}

For a representation $\rho:\g\rightarrow\End(V)$ we define
$\rho^k:\g^{\otimes k} \rightarrow \End(V)$ on simple tensors by
$$ \rho^{k}(g_1\otimes\cdots\otimes g_k)= \rho(g_1) \cdots \rho(g_k)$$
and extend it to the general case by linearity.

\begin{lemma}
\label{lem:highest-weight}
Let $\epsilon$ be a number, let $\rho :
\g \rightarrow \End(V)$ be an irreducible representation of $\g$  and let
$\lambda$ be the corresponding highest weight.

Let $z\in\g^{\otimes k} $ be given; by a small abuse of notation we will denote by
$\rho^k(z)\big|_{V_{\lambda}}: V_{\lambda} \rightarrow V_{\lambda}$
the restriction of $\rho^k(z)$ to $V_{\lambda}$ projected again onto $V_{\lambda}$.
Since the highest-weight space $V_{\lambda}$ in the decomposition \eqref{eq:decomposition2} is one-dimensional, we will identify $\rho^k(z)\big|_{V_{\lambda}}$ with a complex number.
Furthermore, $\epsilon^k \rho^k(z)
\big|_{V_{\lambda}}$ is a polynomial in $\{\Lambda(x):x\in \g\}$ and $\epsilon$,
where 
$\Lambda=\epsilon \lambda$. This polynomial can be written as
$$  \epsilon^k \rho^k(z) \big|_{V_{\lambda}} = \Lambda^{\otimes k}(z)+
(\text{terms of degree at least 1 in $\epsilon$}). $$

\end{lemma}

\begin{proof}
  It is enough to prove that
  \begin{multline*}  \epsilon^k \rho ( g_1 ) \cdots \rho ( g_k )
\big|_{V_{\lambda}} = \\ \begin{cases}
       \Lambda ( g_1 ) \cdots \Lambda ( g_k ) + (\text{terms of degree at least
1 in $\epsilon$}) & \text{if } g_1, \dots, g_k \in
       \h,\\
      (\text{terms of degree at least 1 in $\epsilon$}) & \text{otherwise}
     \end{cases} 
\end{multline*}
  holds true for all tuples $g_1, \dots, g_k \in \g$ such that each $g_i$
  belongs either to $\h$ or to one of the root spaces $\g_{\alpha}$ in the
decomposition \eqref{eq:decomposition}. We 
  use induction over $k$: assume that the lemma holds true for all $k' < k$.

  Let $\pi \in S_k$ be a permutation. Thanks to the commutation relations $xy
  = yx + [ x, y ]$ in the universal enveloping algebra of $\g$ we may write in
the
  universal enveloping algebra:
  \[ g_1 \cdots g_k = g_{\pi ( 1 )} \cdots g_{\pi ( k )} + \text{(summands
     with at most $k - 1$ factors)} \]
  hence
  \begin{multline*}
{\epsilon}^k {\rho}(g_1) {\cdots}
  {\rho}(g_k)  \big|_{V_{\lambda}} = {\epsilon}^k
  {\rho}(g_{{\pi}(1)}) {\cdots} {\rho}(g_{{\pi}(k)})
  \big|_{V_{\lambda}}+\\
  {\epsilon}^k {\rho}(\text{summands with at most $k - 1$ factors})
  \big|_{V_{{\lambda}}}.
\end{multline*}
  The induction hypothesis can be applied to the
  second summand on the right-hand side and it shows that it is of degree at
least $1$ in $\epsilon$, therefore
\begin{multline*}{\epsilon}^k {\rho}(g_1) {\cdots}
  {\rho}(g_k)  \big|_{V_{{\lambda}}} = {\epsilon}^k
  {\rho}(g_{{\pi}(1)}) {\cdots} {\rho}(g_{{\pi}(k)})
  \big|_{V_{{\lambda}}} +\\ (\text{terms of degree at least $1$ in $\epsilon$}). 
\end{multline*}

  It follows that it is enough to consider the case
  when the elements $g_1, \dots, g_k$ are sorted in such a way that for some
  $r, s, t \geq 0$ such that $r + s + t = k$ the initial $r$ elements $g_1,
  \dots, g_r$ belong to root spaces corresponding to the negative roots, the
  next $s$ elements $g_{r + 1}, \dots, g_{r + s}$ belong to $\h$ and the final
  $t$ elements $g_{r + s + 1}, \dots, g_{r + s + t}$ belong to root spaces
  corresponding to the positive roots. A direct calculation shows that
  \[ \epsilon^k \rho ( g_1 ) \cdots \rho ( g_k ) \big|_{V_{
     \lambda}} = \begin{cases}
       \epsilon^k \lambda( g_1 ) \cdots \lambda ( g_k ) & \text{if }
       g_1, \dots, g_k \in \h,\\
       0 & \text{otherwise}
     \end{cases} \]
  and thus the inductive step follows.
\end{proof}

\begin{lemma}
\label{lem:oszacowanie-lambdy}
For a given irreducible representation $\rho^{\lambda}$ of $\g$ we denote
$$M_{\lambda} = - \sum_i \rho^{\lambda}(x_i)^2, $$
where $(x_i)$ denotes an orthogonal basis of $\g$ (regarded as a real vector space).
$M_{\lambda}$ is a multiple of identity hence can be identified with a complex number.

There exists a constant $C$ with a property that
$$ |\lambda|^2 \leq 2 M_{\lambda} + C $$
for any value of $\lambda$.
\end{lemma}
\begin{proof}
Let $(e_i)$ be some linear basis of $\g$ (this time regarded as a complex vector space) and $(f_i)$ be its dual base. Then 
$$M_{\lambda} = - \sum_i \rho^{\lambda}(e_i) \rho^{\lambda}(f_i). $$
It will be convenient for us to take as $(e_i)$ a union of two families: firstly from each root space $\g_{\alpha}$ we
select some non-zero vector (the corresponding dual vector $f_i$ belongs to $\g_{-\alpha}$) and secondly we select some base of $\h$.

Since $M_{\lambda}$ is a multiple of identity Lemma \ref{lem:highest-weight} can be used to 
evaluate it. It is easy to check that there is some element $x\in\h$ (which does not depend on the choice of $\lambda$) with a property that
$$ M_{\lambda} = |\lambda|^2 + \lambda(x).$$
The estimate
$$ 2 \lambda(x)  \geq -|\lambda|^2 - |x|^2 $$
finishes the proof.
\end{proof}

\begin{lemma}
\label{lem:entries_determine_eigenvalues}
Let $\lambda$ be a random vector in $\h$ which is invariant under the action 
of the Weyl group $W$ and let $A$ be the corresponding $G$-invariant random vector in $\g$.
Then each moment $m_k(\lambda)$ is a polynomial function in the moment $m_k(A)$.
\end{lemma}
\begin{proof}
Let us assume for simplicity that $\g$ is semisimple. The Harish-Chandra isomorphism
(see \cite{Knapp2002} for a reference) is an isomorphism between $Z(\bar{\g})$ (the center of the enveloping
algebra $\bar{g}$) and $S(\h)^W$ (the $W$-invariant part of the symmetric algebra $S(\h)$).

The isomorphism of vector spaces $\bar{\g} \cong \bigoplus_{k\geq 0} \Sym^k(\g) $ 
implies that 
\begin{equation} 
\label{eq:HCa}
Z(\bar{\g}) \cong \bigoplus_{k\geq 0} \big[\Sym^k(\g)\big]^G; 
\end{equation}
similarly
\begin{equation}
\label{eq:HCb}
S(\h)^W \cong \bigoplus_{k\geq 0} \big[ \Sym^k(\h) \big]^W.
\end{equation}
The Harish-Chandra isomorphism provides an isomorphism between the right-hand sides of
\eqref{eq:HCa} and \eqref{eq:HCb}. Since it preserves the gradation 
it is the required map.

The general case follows from the fact that the Lie algebra $\g$ can be written 
as a direct sum of its center and a semisimple Lie algebra.
\end{proof}

\subsection{Weak topology and uniqueness}
\label{sec:wtuniqueness}

In view of Theorem \ref{theo:main2} it is natural 
to address the question whether the notion of convergence of moments
considered there
could be replaced by some more probabilistic notion of convergence. 

The notion of convergence in moments is very effective in 
the framework of bounded operators, as it can be seen for example in Voiculescu's
free probability theory \cite{VoiculescuDykemaNica}. 
Indeed, the convergence in moments fully determines the von Neumann algebra generated by
the limiting operator, and in the commutative case the
moments of a bounded random variable determine its distribution by
Stone--Weierstrass
theorem. 

Unfortunately, if the random variables under consideration
are not bounded, then their moments
might not be finite; even if the latter case holds then in general the moments do not
determine the distribution of a random variable. Therefore it would be desirable
to use some more refined description of the joint distribution of random variables.
In the context of classical probability theory such a description is given by
an appropriate probability measure;
unfortunately, it is not clear what would 
be a good notion of distribution of unbounded operators in the non-commutative
case. Some attempts to define weak convergence of joint distribution in the
context of $W\gwia$--probability spaces have been made in very
specific examples (see \cite{MeyerQuantumprobability} and references therein). 
Unfortunately, it is not clear to us how one can adapt these definitions 
in our setting. A promising approach to this problem via Gromov-Hausdorff
distance was presented by Rieffel \cite{Rieffel2004}, however 
for the moment it is
not clear if this method can be succesfully used for our purposes.

To summarize the above discussion: we have no candidate for some kind of convergence
which would replace the convergence in moments in point (\ref{convergence:B}) of Theorem \ref{theo:main2}
since we deal here with a joint distribution of a family of non-commuting random variables.

Nevertheless, in point \eqref{convergence:A} of Theorem \ref{theo:main2} we deal with classical random 
variables therefore it makes sense to consider convergence in some other sense, such as weak
convergence of probability measures. 
Unfortunately, in general there is no connection between convergence  in moments and weak convergence of probability measures. 
In particular, in order for convergence of measures in moments to imply 
their weak convergence we need to assume, for example, that the limit measure is
uniquely determined by its moments.

It seems to be hard to incorporate
the weak convergence of probability measures to condition \eqref{convergence:A} and
preserve the equivalence of conditions \eqref{convergence:A} and \eqref{convergence:B}. Therefore our strategy 
will be to keep Theorem \ref{theo:main2} unchanged and in the study of its 
applications to pay attention to the weak convergence of probability measures
(e.g. Theorem \ref{theo:restriction}, item (2) and Theorem \ref{theo:tensor}, item (2)).

\section{Applications of the main theorem}
\label{sec:applications}

\subsection{Restriction of representations and tensor product of representations}
\label{sec:restriction-tensorproduct}

In this section we investigate a few remarkable consequences of 
Theorem \ref{theo:main2}.

\begin{theorem}
\label{theo:restriction}
Let $G\subset G'$ be Lie groups and $\g\subset \g'$ be the corresponding Lie
algebras, let $(\epsilon_n)$ be a sequence of real numbers which converges to
zero. Let $(\rho'_n)$ be a sequence of representations of $\g'$;
by $\lambda'_n$ and $\lambda_n$ we denote the random highest weight 
corresponding to representations $\rho'_n$ and $\rho_n'|_{\g}$, respectively.

\renewcommand{\labelenumi}{(\arabic{enumi})}
\begin{enumerate}
\item \label{cond-moms} 

Assume that $\epsilon_n \lambda'_n$ converges in moments to the distribution
of a $G'$-invariant random vector $A$ with values in $\g'$. Then the sequence 
$\epsilon_n \lambda_n$ converges in moments towards the $G$-invariant random
vector $\Pi_{\g}(A)$ where $\Pi_{\g}:\g'\rightarrow\g$ is the orthogonal
projection.

\item \label{cond-weak}
Assume that $\epsilon_n \lambda'_n$ converges weakly to the distribution of a $G'$-invariant random vector $A$ with values in $\g'$.
Then the sequence 
$\epsilon_n \lambda_n$ converges weakly towards the $G$-invariant random
vector $\Pi_{\g}(A)$.

\end{enumerate}

\end{theorem}

\begin{proof}
Notice that the non-commutative random vector $\rho'_n|_{\g}$ is a
projection of the non-commutative random vector $\rho'_n$ onto $\g$.
It follows that the random matrix which is the entrywise limit of $\epsilon_n
\rho'_n |_{\g}$ is  a projection of $A$ to $\g$. Theorem \ref{theo:main2} can be
applied twice: for the sequence $(\rho'_n)$ and for the sequence $(\rho'_n|_{\g})$
which finishes the proof of part (\ref{cond-moms}).

In order to prove (\ref{cond-weak}) it is enough to show that for every
$\varepsilon>0$ and every subsequence $\epsilon_{k(n)} \lambda_{k(n)}$
we can chose a subsequence $\epsilon_{k(l(n))} \lambda_{k(l(n))}$ which converges
weakly to some limit distribution with a property that its variation distance
from the distribution of $\Pi_{\g}(A)$ is smaller than $\varepsilon$.

Let $\varepsilon>0$ and a subsequence $\epsilon_{k(n)} \lambda_{k(n)}$ be fixed. 
For simplicity, in the following instead of $\epsilon_{k(n)} \lambda_{k(n)}$ we
shall consider just the sequence $\epsilon_{n} \lambda_{n}$.
We can find a sequence of representations
$(\tilde{\rho}'_n)$ for which corresponding rescaled random highest weights
$\epsilon_n \tilde{\lambda}'_n$ have a common compact support and 
the total variation distance between the distribution of $\epsilon_n \lambda'_n$ and 
$\epsilon_n \tilde{\lambda}'_n$ is smaller than $\varepsilon$  
(such a sequence $\tilde{\rho}'_n$ can be constructed by truncating the distribution
of $\epsilon_n \tilde{\lambda}'_n$ to some sufficiently big compact set).

By compactness argument we can select a subsequence $\epsilon_{l(n)} \tilde{\lambda}'_{l(n)}$
which converges weakly (hence in moments) to some limit; let $\tilde{A}$ be a
$G$-invariant random vector in $\g\gwia$ with this distribution of eigenvalues.
The total variation distance between the distribution of $\tilde{A}$ and $A$ is bounded
by $\varepsilon$.

The first part of the theorem can be applied to the subsequence of representations
$\epsilon_{l(n)} \tilde{\rho}_{l(n)}$; it follows that the highest weights
$\epsilon_{l(n)} \tilde{\lambda}_{l(n)}$ (corresponding to the restrictions of
$\tilde{\rho}_{l(n)}$ to $\g$) converge in moments to the
distribution of the eigenvalues of $\Pi_{\g}(\tilde{A})$. Since the random weights
have a common compact support the convergence holds also in the weak sense.
We finish the proof by observing that the total variation distance between
the distribution of $\Pi_{\g}(\tilde{A})$ and $\Pi_{\g}(A)$ is smaller than $\varepsilon$.
\end{proof}

In the case of the inclusion of the groups $U(d)\subseteq U(d')$ for $d<d'$
the above theorem takes the following concrete form:

\begin{corollary}
\label{coro:restriction}
Let $d<d'$ be positive integers and let $(\epsilon_n)$ be a sequence of real
numbers which converges to zero. Let $A=(A_{ij})_{1\leq i,j\leq d'}$ be a hermitian $U(d)$-invariant random matrix and let $\rho_n$ be a sequence of representations
of\/ $U(d')$ with a property that the distribution of of $\epsilon_n \lambda'_n$ converges to the joint distribution of eigenvalues of $A$, where $\lambda'_n\in \Z_{d'}$ is a random weight associated to $\rho_n$.

Then the distribution of of $\epsilon_n \lambda_n$ converges to the joint distribution of eigenvalues of the corner $(A_{ij})_{1\leq i,j\leq d}$, where $\lambda_n\in \Z_{d}$ is a random weight associated to the restriction $\epsilon_n  \rho_{n}|_{\uu_{d}} $.
\end{corollary}

Observe that in the above result we used Remark \ref{re:Un} in order to work with hermitian random matrices. Similar concrete interpretation in the case of unitary groups is also possible for the other results presented in this section.

We leave to the Reader to investigate the other such simple inclusions: of
orthogonal groups $O(d)\subset O(d')$ and of symplectic groups $Sp(d)\subset
Sp(d')$.

\begin{theorem}
\label{theo:tensor}
Let $(\epsilon_n)$ be a sequence of real numbers which converges to zero,
let $(\rho^{(1)}_n)$, $(\rho^{(2)}_n)$ be two sequences of representations of
$\g$ and let $\lambda^{(1)}_n$ and $\lambda^{(2)}_n$ be the corresponding sequences
of random highest weights. Furthermore, let $\lambda_n$ be the sequence of random
highest weights corresponding to the tensor products
$\rho^{(1)}_n\otimes \rho^{(2)}_n$.
Let $A^{(1)}$, $A^{(2)}$ be independent $G$-invariant random matrices in
$\g\gwia$.

\begin{enumerate}
 \item \label{tensor:moms}
If for each $i\in\{1,2\}$ the sequence $\epsilon_n \lambda^{(i)}_n$
converges in moments to the distribution of eigenvalues of $A^{(i)}$
then
$ \epsilon_n \lambda_n $ converges in moments to the distribution of
eigenvalues of $A^{(1)}+A^{(2)}$.

 \item \label{tensor:weak}
If for each $i\in\{1,2\}$ the sequence $\epsilon_n \lambda^{(i)}_n$
converges weakly to the distribution of eigenvalues of $A^{(i)}$
then
$ \epsilon_n \lambda_n $ converges weakly to the distribution of
eigenvalues of $A^{(1)}+A^{(2)}$.

\end{enumerate}

\end{theorem}
\begin{proof}
Let $\rho_n^{(3)}:=\rho_n^{(1)} \otimes \rho_n^{(2)} $. Then
\begin{equation}
\label{eq:iloczyn-tensorowy}
 \rho_n^{(3)}(x)= \rho_n^{(1)}(x) \otimes 1 +1 \otimes  \rho_n^{(2)}(x).
\end{equation}
It follows that $\epsilon_n \rho_n^{(3)}$ viewed as a non-commutative random
vector is a sum of two non-commutative random vectors:
$\big(\epsilon_n \rho_n^{(1)}(x)\big) \otimes 1$ and $1\otimes
\big( \epsilon_n \rho_n^{(2)}(x) \big)$. Theorem \ref{theo:main2} implies that
the first summand converges in moments to $A^{(1)}$ and the second
summand converges in moments to $A^{(2)}$. The coordinates of the first vector (viewed
as non-commutative random variables) commute with the coordinates of the second
vector; since we consider them with respect to a state (i.e.\ the
normalized trace) which is a tensor product of the original
states (i.e.\ normalized traces), it follows that their sum converges to the sum of the independent random matrices which finishes the proof of part (\ref{tensor:moms}).

Part (\ref{tensor:weak}) follows from part (\ref{tensor:moms}) in a similar
way as in Theorem \ref{theo:restriction}.
\end{proof}

\subsection{Central limit theorem}

The following theorem was already proved by Kuperberg \cite{Kuperberg2002} in a
slightly different setting.
We include this theorem here because we believe that the proof in this
framework 
(using Theorem \ref{theo:main2}) 
is new.

\begin{theorem}
\label{theo:clt}
Let $G$ be a Lie group and let $\rho$ be its representation.  
We define $c\in\g\gwia$ given by $c(x)=\tr \rho(x)$.

Let $\lambda_n$ be the random highest weight corresponding to the representation
$\rho^{\otimes n}$.
Then the sequence 
\begin{equation}
\label{eq:clt-for-highest-weight}
\frac{1}{\sqrt{n}} \big[ \lambda_n - c n\big]
\end{equation}
converges (both in moments and weakly) to the distribution of eigenvalues of a
certain 
centered Gaussian random matrix in $\g\gwia$.

The covariance of the above Gaussian random matrix is given by
$$ \int_{\g\gwia} \lambda(x)\ \lambda(y) d\mu(\lambda) = \tr \rho(x) \rho(y) $$
for any $x,y\in\g$. 
\end{theorem}
\begin{proof}
Similarly as in Eq.\ \eqref{eq:iloczyn-tensorowy} we have that
$$ \rho^{\otimes n}(x)= \rho(x) \otimes 1 \otimes \cdots +1 \otimes  \rho(x)
\otimes 1 \otimes \cdots + \cdots. $$
regarded as a non-commutative random variable is a sum of $n$ commuting
summands therefore the non-commutative central limit theorem of Giri and
von Waldenfels \cite{GiriWaldenfels1978} can be applied. It follows that the
distribution of the non-commutative random vector $\frac{1}{\sqrt{n}} \big[
\rho^{\otimes n} - c n\big] $ converges in moments to the distribution of
a certain centered Gaussian random variable $X$ which takes values in $\g\gwia$.
We apply Theorem \ref{theo:main2}; it follows that the distribution of the
random weight \eqref{eq:clt-for-highest-weight} converges in moments to the
distribution of the eigenvalues of $X$.

Let $\|\cdot\|$ be the norm on $\g\gwia$ associated to any $G$-invariant scalar
product. The distribution of $X$ is multidimensional Gaussian therefore
there exists a constant $d>0$ such that the tail estimate 
\begin{equation}
\label{eq:tail}
P\left(\|X\| > t\right) < e^{-d t^2} 
\end{equation}
holds true for sufficienly big values of $t$.  For any element of $\g\gwia$ its
norm is equal to the norm of the element of $\h\gwia$ corresponding to its
eigenvalues; it follows that the estimate \eqref{eq:tail} remains true if the
random variable $X$ is replaced by its eigenvalues. The estimate
\eqref{eq:tail} shows therefore that the even moments of the eigenvalues
distribution of $X$ are dominated by the even moments of a Gaussian
distribution.

This implies that $\E(e^{\|X\|})<\infty$ and by Corollary 2.2 and Theorem 2.1
of \cite{DvurevcenskijLahtiYlinen2002} it follows that
 the distribution of eigenvalues of $X$ is
uniquely determined by its moments. It follows by a standard compactness
argument that the distribution of the
random weights \eqref{eq:clt-for-highest-weight} converges weakly to the
distribution of the eigenvalues of $X$. 
\end{proof}

Often we have some additional information about the considered Lie group $G$ which restricts the number of $G$--invariant Gaussian measures on $\g\gwia$. In the case of $G=U(d)$ it is convenient to consider a centered hermitian Gaussian random matrix $g=(g_{ij})_{1\leq i,j\leq d}$ defined by the convariance
$$ \E g_{ij} g_{kl}=0, \qquad \E g_{ij} \overline{g_{kl}}=\delta_{il} \delta_{jk}; $$
this kind of random matrix (and the corresponding measure on the space of $d\times d$ hermitian matrices) is called \emph{Gaussian Unitary Ensemble} ($\text{GUE}$) and plays an important role in random matrix theory.
For any $v\geq 0$ we define $\text{GUE}_v$ as the distribution of a random matrix
$$ g - \tr g + x $$
where $x$ is an independent centered Gaussian variable with the variance $v$. 

Under the isomorphism from Remark \ref{re:Un}, $\text{GUE}_v$ becomes a measure on the Lie algebra $u(d)$
and it is not very difficult to check that (except for degenerate cases) every $U(d)$--invariant Gaussian measure (up to dilation by some number) is of this form.
In particular, we get the following result.

\begin{corollary}
Let $\rho$ be a representation of $U(d)$.
There exist constants $c_1,c_2$ with a property that if $\lambda_n=(\lambda_{n,1},\dots,\lambda_{n,d})$ is a random weight associated to $\rho^{\otimes n}$ then the joint distribution of the components of the vector
$$\frac{c_1}{\sqrt{n}} ( \lambda_n - n c_2) $$
converges to the joint distribution of the eigenvalues of the $\text{GUE}_v$ random matrix.
\end{corollary}

\subsection{Toy example: representations of $SO(3)$ and $SU(2)$}

The above results in the simplest non-trivial case of $G=SO(3)$ should
not be very surprising from a viewpoint of quantum mechanics. Each
quantum-mechanical system in three-dimensional space can be viewed as a
(possibly reducible)
representation of $SO(3)$ (or its universal cover $\Spin(3)=SU(2)$) on some
Hilbert space $V$. The
irreducible components of this representation have a nice physical
interpretation as physical states with a well-defined length $|J|$ of the
angular momentum. For simplicity,  we assume that
$V$ itself is irreducible hence $V$ is finite-dimensional.
The information about the state of the physical system is
encoded by a state $\phi$  on the algebra generated
by observables. We are interested in the situation when the physical state of
the system is $SO(3)$-invariant; it follows that $\phi$ is the normalized trace
on $\End(V)$. 

The physicist's question about the distribution of a component $J_z$ of the
angular momentum can be reformulated in the language of mathematics as question
about the decomposition into irreducible components of the restriction
$V\downarrow^{SO(3)}_{SO(2)}$ of the representation $V$ to a subgroup $SO(2)$
(or, more generally, restriction of the representation $V$ of $\Spin(3)=SU(2)$
to its subgroup $\Spin(2)=U(1)$), 
namely it is the uniform measure on the set
of integers (or half-integers) 
\begin{equation} \label{eq:jz}
\big\{ -|J|,-|J|+1,\dots,|J|-1,|J| \big\}
\end{equation}
(for simplicity we use the system of units in which the Planck's constant
$\hbar=1$).

On the other side it is well-known that when the size of our system
becomes macroscopic then quantum mechanics may be approximated by the classical
mechanics, where the angular momentum $\vec{J}=(J_x,J_y,J_z)$ is just a
usual vector consisting of numbers.  Our assumptions on irreducibility 
and $SO(3)$-invariance imply that in the classical limit $\vec{J}$ is a random
vector with a uniform distribution on a sphere of a fixed length $|J|$.
The physicists question about the distribution of a component $J_z$ of the
angular momentum in this context is answered by the theorem of Archimedes,
namely it is the uniform distribution on the interval $[-|J|,|J|]$. 
In the limit $|J|\to\infty$, after appropriate rescaling,
the uniform measure on the set \eqref{eq:jz} indeed converges to the uniform
measure on the set $[-|J|,|J|]$ hence the answer given by classical mechanics
is indeed the limit of the answer given by quantum mechanics.

It is also a consequence of this article: indeed we view the
representation $\rho$ of the Lie algebra $so(3)$ as an element of $\big( so(3)
\big)\gwia \otimes \End(V)$; in our current case this takes a concrete form of
a matrix
\begin{equation} \label{eq:macierz-j}
J = \begin{bmatrix}
0     & J_z  & -J_y  \\
-J_z  & 0    & J_x   \\
J_y   & -J_x & 0   
\end{bmatrix}
\end{equation}
(the form of this matrix depends on particular choice of the identification of
$\big( so(3)\big)\gwia$ with certain $3\times 3$ matrices).
Theorem \ref{theo:main2} states that asymptotically the matrix
\eqref{eq:macierz-j} behaves like a random $3\times 3$ antisymmetric matrix
with eigenvalues $0, |J|, -|J|$.

We leave the analysis of the central limit theorem in this case to the Reader.

\section{Acknowledgments} \label{sec:Acknowl}The research of P.\'S.\ was
supported by State Committee for Scientific Research (KBN) grant \text{2 P03A
007 23}, RTN network: QP-Applications contract No.~HPRN-CT-2002-00279, and
KBN-DAAD project 36/2003/2004.

The research of B.C.\ was partly supported by a RIMS fellowship and by CNRS.

Acknowledgment: the authors thank the referee for
his feedback and suggestions of improvement to the paper.

\bibliographystyle{alpha}

\bibliography{biblio}

\begin{thebibliography}{KSW96b}

\bibitem[Bia95]{Biane95}
Philippe Biane.
\newblock Representations of unitary groups and free convolution.
\newblock {\em Publ. Res. Inst. Math. Sci.}, 31(1):63--79, 1995.

\bibitem[Bia98]{Biane1998}
Philippe Biane.
\newblock Representations of symmetric groups and free probability.
\newblock {\em Adv. Math.}, 138(1):126--181, 1998.

\bibitem[BtD95]{Brocker-tomDieck}
Theodor Br{\"o}cker and Tammo tom Dieck.
\newblock {\em Representations of compact {L}ie groups}, volume~98 of {\em
  Graduate Texts in Mathematics}.
\newblock Springer-Verlag, New York, 1995.

\bibitem[C{\'S}07]{CollinsSniady2007}
Beno{\^{\i}}t Collins and Piotr {\'S}niady.
\newblock Representations of {L}ie groups, random matrices and free
  probability.
\newblock In preparation, 2007.

\bibitem[DLY02]{DvurevcenskijLahtiYlinen2002}
Anatolij Dvure{\v{c}}enskij, Pekka Lahti, and Kari Ylinen.
\newblock The uniqueness question in the multidimensional moment problem with
  applications to phase space observables.
\newblock {\em Rep. Math. Phys.}, 50(1):55--68, 2002.

\bibitem[FH91]{FultonHarris}
William Fulton and Joe Harris.
\newblock {\em Representation theory}, volume 129 of {\em Graduate Texts in
  Mathematics}.
\newblock Springer-Verlag, New York, 1991.

\bibitem[GvW78]{GiriWaldenfels1978}
N.~Giri and W.~von Waldenfels.
\newblock An algebraic version of the central limit theorem.
\newblock {\em Z. Wahrscheinlichkeitstheorie und Verw. Gebiete},
  42(2):129--134, 1978.

\bibitem[GW98]{GoodmanWallach}
Roe Goodman and Nolan~R. Wallach.
\newblock {\em Representations and invariants of the classical groups},
  volume~68 of {\em Encyclopedia of Mathematics and its Applications}.
\newblock Cambridge University Press, Cambridge, 1998.

\bibitem[Hec82]{Heckman1982}
G.~J. Heckman.
\newblock Projections of orbits and asymptotic behavior of multiplicities for
  compact connected {L}ie groups.
\newblock {\em Invent. Math.}, 67(2):333--356, 1982.

\bibitem[Joh01]{Johansson01}
Kurt Johansson.
\newblock Discrete orthogonal polynomial ensembles and the {P}lancherel
  measure.
\newblock {\em Ann. of Math. (2)}, 153(1):259--296, 2001.

\bibitem[Kna02]{Knapp2002}
Anthony~W. Knapp.
\newblock {\em Lie groups beyond an introduction}, volume 140 of {\em Progress
  in Mathematics}.
\newblock Birkh\"auser Boston Inc., Boston, MA, second edition, 2002.

\bibitem[KSW96a]{KazakovStaudacherWynter1}
Vladimir~A. Kazakov, Matthias Staudacher, and Thomas Wynter.
\newblock Almost flat planar diagrams.
\newblock {\em Comm. Math. Phys.}, 179(1):235--256, 1996.

\bibitem[KSW96b]{KazakovStaudacherWynter2}
Vladimir~A. Kazakov, Matthias Staudacher, and Thomas Wynter.
\newblock Character expansion methods for matrix models of dually weighted
  graphs.
\newblock {\em Comm. Math. Phys.}, 177(2):451--468, 1996.

\bibitem[Kup02]{Kuperberg2002}
Greg Kuperberg.
\newblock Random words, quantum statistics, central limits, random matrices.
\newblock {\em Methods Appl. Anal.}, 9(1):99--118, 2002.

\bibitem[Kup05]{Kuperberg2005}
Greg Kuperberg.
\newblock A tracial quantum central limit theorem.
\newblock {\em Trans. Amer. Math. Soc.}, 357(2):459--471 (electronic), 2005.

\bibitem[Lit95]{Littelmann95}
Peter Littelmann.
\newblock Paths and root operators in representation theory.
\newblock {\em Ann. of Math. (2)}, 142(3):499--525, 1995.

\bibitem[Meh91]{Mehta}
Madan~Lal Mehta.
\newblock {\em Random matrices}.
\newblock Academic Press Inc., Boston, MA, second edition, 1991.

\bibitem[Mey93]{MeyerQuantumprobability}
Paul-Andr{\'e} Meyer.
\newblock {\em Quantum probability for probabilists}, volume 1538 of {\em
  Lecture Notes in Mathematics}.
\newblock Springer-Verlag, Berlin, 1993.

\bibitem[Rie04]{Rieffel2004}
Marc~A. Rieffel.
\newblock {\em Gromov-{H}ausdorff distance for quantum metric spaces. {M}atrix
  algebras converge to the sphere for quantum {G}romov-{H}ausdorff distance}.
\newblock American Mathematical Society, Providence, RI, 2004.
\newblock Mem. Amer. Math. Soc. {\bf 168} (2004), no. 796.

\bibitem[{\'S}ni06]{Sniady2005GaussuanFluctuationsofYoungdiagrams}
Piotr {\'S}niady.
\newblock {G}aussian fluctuations of characters of symmetric groups and of
  {Y}oung diagrams.
\newblock {\em Probab. Theory Related Fields}, 136(2):263--297, 2006.

\bibitem[{\'S}S07]{SniadySpeicher2007permutationmatrix}
Piotr {\'S}niady and Roland Speicher.
\newblock Permutationally invariant random matrices.
\newblock in preparation, 2007.

\bibitem[VDN92]{VoiculescuDykemaNica}
D.~V. Voiculescu, K.~J. Dykema, and A.~Nica.
\newblock {\em Free random variables}.
\newblock American Mathematical Society, Providence, RI, 1992.

\bibitem[Wor87]{Woronowicz87Compact-matrix-pseudogroups}
S.~L. Woronowicz.
\newblock Compact matrix pseudogroups.
\newblock {\em Comm. Math. Phys.}, 111(4):613--665, 1987.

\end{thebibliography}

\end{document}